%% file: main.tex
\documentclass{article}

\input{packages.tex}

\title{PolyTO: Structural Topology Optimization using Convex Polygons}

\author{
Aaditya Chandrasekhar \\
Department of Mechanical Engineering\\
University of Wisconsin-Madison\\
Madison, WI 53706 \\
  \texttt{cs.aaditya@gmail.com} \\
}

\begin{document}
\newcommand\barbelow[1]{\stackunder[1.2pt]{$#1$}{\rule{.8ex}{.075ex}}}
\maketitle
\begin{abstract}

In this paper, we propose a topology optimization (TO) framework where the design is parameterized by a set of convex polygons. Extending feature mapping methods in TO, the representation allows for direct extraction of the geometry. In addition, the method allows one to impose geometric constraints such as feature size control directly on the polygons that are otherwise difficult to impose in density or level set based approaches. The use of polygons provides for more more varied shapes than simpler primitives like bars, plates, or circles. The polygons are defined as the feasible set of a collection of halfspaces.  Varying the halfspace's parameters allows for us to obtain diverse configurations of the polygons. Furthermore, the halfspaces are differentiably mapped onto a background mesh to allow for analysis and gradient driven optimization. The proposed framework is illustrated through numerous examples of 2D structural compliance minimization TO. Some of the key limitations and future research are also summarized.

 \begin{figure*}[htbp]
 	\begin{center}
		\includegraphics[scale=1.2,trim={0 20 0 0},clip]{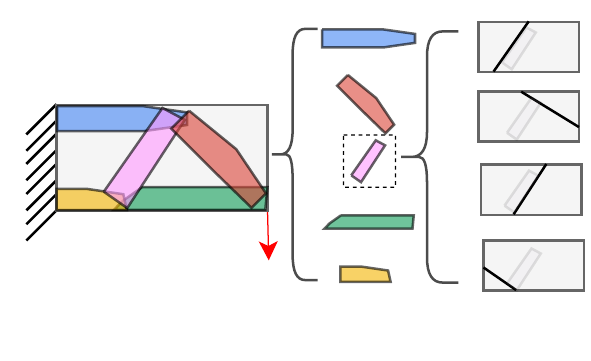}
 		\caption{An optimized design constituted by a set of convex polygons. Each polygon in turn is defined by the feasible set of a collection of half-spaces.}
 		\label{fig:abstract_polyto}
	\end{center}
 \end{figure*}
 
\end{abstract}


\section{Introduction}
\label{sec:intro}
\input{intro.tex}

\section{Proposed Method}
\label{sec:method}
\input{method.tex}

\section{Numerical Results}
\label{sec:results}
\input{results.tex}

\section{Conclusion}
\label{sec:conclusion}

\input{conclusion.tex}

\section*{Compliance with ethical standards}
\label{sec:ethics}
The authors declare that they have no conflict of interest.

\section*{Replication of Results}
\label{sec:replication}
The Python code pertinent to this paper can be found at \href{https://github.com/aadityacs/PolyTO}{github.com/aadityacs/PolyTO}

\bibliographystyle{unsrt}  
\bibliography{references} 
\end{document}

%% file: packages.tex
\usepackage{arxiv}

\usepackage[utf8]{inputenc} 
\usepackage[T1]{fontenc}    
\usepackage{hyperref}       
\usepackage{url}            
\usepackage{algpseudocode}
\usepackage{amsmath}
\usepackage{booktabs}       
\usepackage{amsfonts}       
\usepackage{nicefrac}       
\usepackage{caption}
\usepackage{subcaption}
\usepackage{microtype}      
\usepackage{lipsum}
\usepackage{graphicx}
\usepackage{bm}
\usepackage{algorithm}
\usepackage{cleveref}
\usepackage{xcolor}
\usepackage{stackengine}
\graphicspath{ {./images/} }

%% file: intro.tex
Topology optimization (TO) is a strategy for optimally distributing material within a design domain \cite{sigmund2013TOReview} to optimize a desired objective, while meeting one or more constraints. TO has matured to the point where it is used in commercial settings \cite{meng2020TO_commercial_example} to explore novel and optimal designs. While various TO methods have been proposed; density \cite{bendsoe2003topology} and level set based methods \cite{wang2003levelsetTO} methods have been widely adopted. These methods produce complex, organic designs across a wide range of applications. However, the designs produced are often not directly manufacturable, and require either significant post processing \cite{subedi2020GeomPostProcessing} or imposing corresponding fabrication constraints \cite{vatanabe2016ManufCons} during optimization. Post-processing TO designs may lead to significant degradation of its performance and in cases such as design for novel applications such as compliant mechanisms \cite{zhu2020dCompliantMechanismTO, zhu2013TOCompliantMechanism} may render the design entirely non-functional. This has restricted TO to be seen more as a 'conceptual tool', limiting its full potential in engineering workflows. While imposing fabrication constraints can improve manufacturability, the development of these constraints are research intensive and is the focus of much ongoing research \cite{Ian2022inverse, hammond2021photonic, sutradhar2017incorporating}. Finally, while additive manufacturing (AM) promises to directly produce TO free-form designs, despite requiring additional AM specific fabrication constraints \cite{thompson2016AMConstraints}, many issues remain unresolved \cite{liu2018AMinTOReview} and wider adoption of AM is yet to be realized.

In parallel, TO techniques have been developed recognizing that a significant portion of components manufactured today is done via the composition of geometrically simpler stock material. Broadly termed as feature mapping methods \cite{wein2020reviewFeatureMapping}, it includes TO techniques such as the method of Moving Morphable Voids \cite{zhang2017MMV}, the moving morphable components method \cite{guo2014MMC}, the adaptive bubble method \cite{cai2020adaptiveBubbleMethod}, and the geometric projection \cite{norato2015gpto} method. A distinguishing feature of these techniques is that the design is parameterized by a high-level geometric description with features that are mapped onto a mesh for analysis. This is in contrast to density and level set based method where a the design is parameterized verbosely by pixels/voxels. For a detailed review of these methods, we refer readers to \cite{wein2020reviewFeatureMapping}.

Parameterizing the design via high-level geometric description simplifies the enforcement of rules concerning the presence of and dimensions of features. In addition, various production considerations such as dimensioning and tolerancing, design variability quantification and control require a high-level description of the geometry. Finally, they can be easily converted into polygonal meshes via a simple duality transform \cite{deng2020cvxnet}. This is in contrast to pixel based representations, where a computationally intensive iso-surfacing operation such as marching cubes \cite{lorensen1987marchingCubes} needs to be performed to extract contours.

\subsection{Contributions}
\label{sec:intro_contributions}

In this paper, we present an extension to feature mapping methods where our design is parameterized using convex polygons. Compared to simpler primitives such as bars and plates \cite{zhang2016GPTO_plates}, the parameterization allows for more design flexibility while still being recoverable as geometric primitives. The intrinsic guarantee of the convexity of the primitives further promotes manufacturability. In contrast to \cite{norato2018TO_supershapes}, we showcase that feature sizes constraints can be derived and imposed easily using the proposed parameterization. Further, the signed distance field is defined implicitly and is straightforward to compute.
 
 The remainder of the paper is organized as follows. In Section~\ref{sec:method_designParameterization} we introduce the proposed design parameterization. Section~\ref{sec:method_geometryProjection} describes the computation of the signed distance field (SDF) followed by projection of the SDF onto a density field for the analysis. Section~\ref{sec:method_opt} and Section~\ref{sec:sensitivity_analysis} describe the formulation of the optimization and sensitivity analysis respectively. We summarize the complete algortithm in Section~\ref{sec:method_algo}. Several parametric studies that demonstrate the proposed method are presented in Section~\ref{sec:results}.  Finally,  we discuss some of the limitations and directions for future research in Section~\ref{sec:conclusion}.

%% file: method.tex
\subsection{Overview}
\label{sec:method_overview}

We will assume that a design domain with loads and restraints have been prescribed. Further, we assume that the hyper-parameters: number of polygons $(K)$ and the maximum number of sides of each polygon $(S)$ have been specified. The objective then, is to find a configuration of the polygons that results in an optimal design. For simplicity, we restrict our discussion to the optimization of minimal compliance designs under a volume constraint.

\subsection{Design Parameterization}
\label{sec:method_designParameterization}

Let us suppose our design is to be defined using $K (\geq 1)$ convex polygons. Further, let us suppose each polygon can have a maximum of $S (\geq 3)$ sides. Observe that a convex polygon can be defined as the feasible set of a collection of half-spaces \cite{boyd2004convex}. Recall that while the feasible set defined by a collection of half-spaces is guaranteed to be convex, there is no guarantee of it being bounded or non-empty. To prevent unboundedness and non-empty region, we present a conservative parameterization rather than parameterizing our geometry by the half-space's normal and distance from the origin. Specifically, each polygon $p_i$ is parameterized as:

\begin{equation}
    p_i = \{c_x^{(i)}, c_y^{(i)}, \alpha^{(i)}, d_1^{(i)}, \ldots, d_S^{(i)} \} \quad , \; i = 1,\ldots, K
    \label{eq:param_of_polygon}
\end{equation}

\begin{figure}[h]
 	\begin{center}
		\includegraphics[scale=0.5,trim={0 0 0 0},clip]{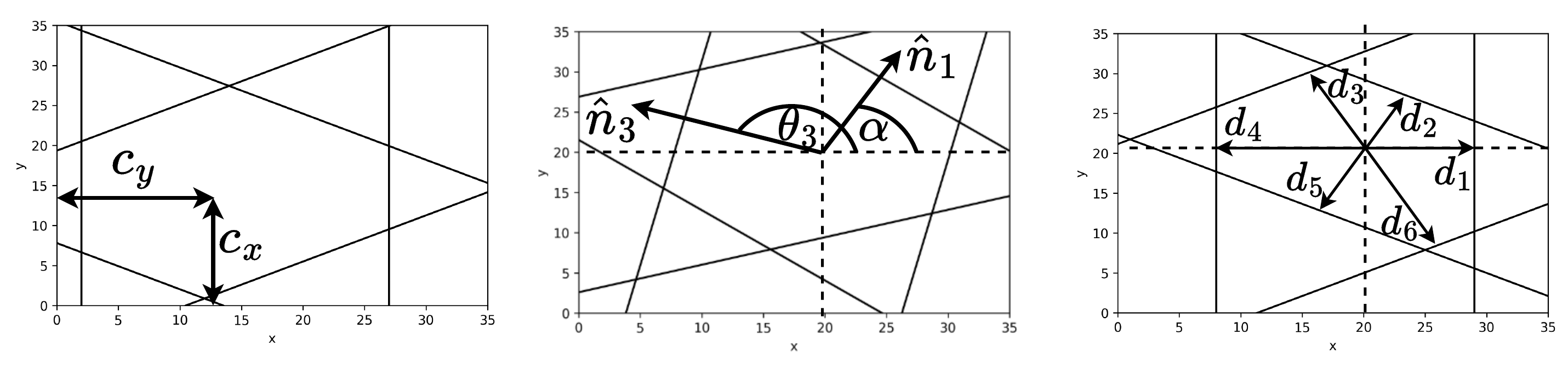}
 		\caption{The (hexagonal) polygon is parameterized by the center coordinate $(c_x, c_y)$, an angular offset $(\alpha)$, and the distances of the half-spaces from the center $(d_1, \ldots, d_6).$}
 		\label{fig:polygon_params}
	\end{center}
 \end{figure}

Where $c_x^{(i)}, c_y^{(i)}$ are respectively the $x$ and $y$ center coordinates of the $i^{th}$ polygon, $d_j^{(i)} \in \mathrm{R}^{+}$ is the distance of the $j^{th}$ half-space of the $i^{th}$ polygon from its center. We define $\alpha_i$ as the overall angular offset of the $i^{th}$ polygon  with the orientation of the $j^{th}$ half-space defined as $\theta_j^{(i)} = \alpha_i + \frac{2\pi(j-1)}{S} $ (Fig.~\ref{fig:polygon_params}). Observe that while this definition of orientation limits the design space, it ensures that the polygon is bounded and non-empty.

\subsection{Geometry Projection}
\label{sec:method_geometryProjection}

The central idea of geometry projection is to map the design as described by the polygon's parameters (Section~\ref{sec:method_designParameterization}) onto a density field defined over a mesh  \cite{deng2020cvxnet}. The density function $\rho(x),x\in\mathbb{R}^2$ is defined such that $\rho(x) = 0$ defines the void regions and $\rho(x) = 1$ the solid regions. A point is to be denoted as solid if it is within the feasible region of at-least one of the polygons. This mapping allows for subsequent structural analysis over the mesh that is required for the optimization. Further the mapping is a differentiable so that sensitivities of the objectives and constraints are defined, enabling gradient-based optimization. 

We begin by computing the signed distance field (SDF) $(\hat{\phi}^{(i)}_j(\cdot))$ of the $j^{th}$ half-spaces of the $i^{th}$ polygon at a point $(x,y)$ in the domain as (Eq.~\ref{eq:sdf_hyperplane}): 

 \begin{equation}
    \hat{\phi}_j^{(i)}(x, y; c_x^{(i)}, c_y^{(i)}, \theta_j^{(i)}, d_j^{(i)}) = (x - c_x^{(i)})\cos(\theta_j^{(i)}) + (y - c_y^{(i)})\sin(\theta_j^{(i)}) - d_j^{(i)}
    \label{eq:sdf_hyperplane}
\end{equation}

\begin{figure}[h]
 	\begin{center}
		\includegraphics[scale=0.6,trim={0 0 0 0},clip]{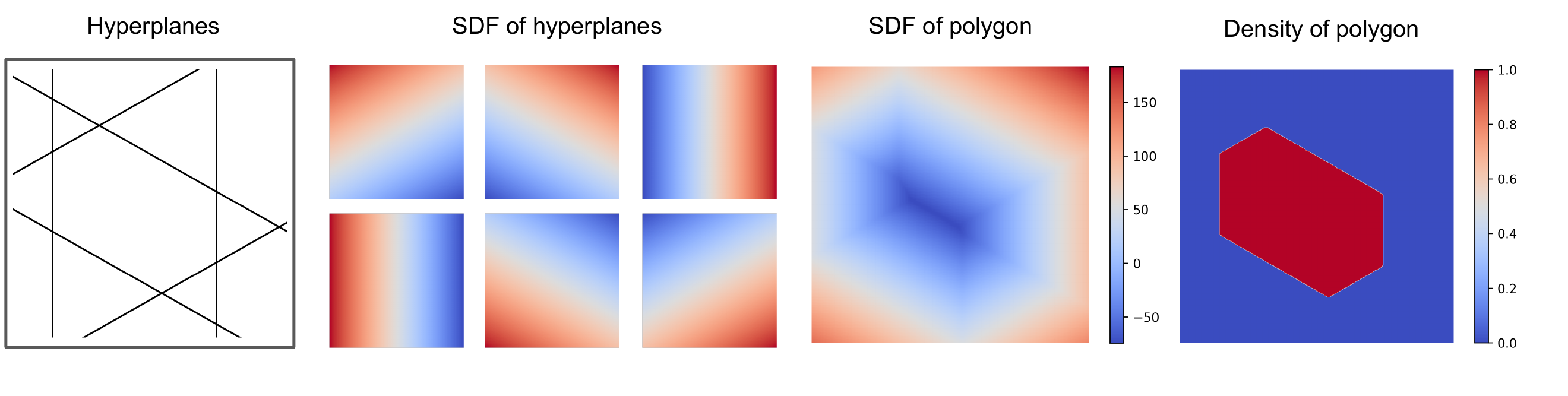}
 		\caption{From half-spaces to density}
 		\label{fig:hyperplane_to_density}
	\end{center}
 \end{figure}

 With the SDF of the half-spaces defined, the $\textit{approximate}$ SDF of the convex polygon can be computed by taking the intersection (max) of the SDF of the planes. To facilitate gradient based learning, instead of maximum, we use the smooth maximum LogSumExp (LSE) function. This results in the approximate SDF of the $i^{th}$ polygon $(\phi(\cdot))$ as (Eq.~\ref{eq:sdf_polygon}):

\begin{equation}
    \phi^{(i)}(x, y; c_x^{(i)}, c_y^{(i)}, \alpha^{(i)}, d_1^{(i)}, \ldots, d_S^{(i)}) = \text{LSE}(\hat{\phi}_1^{(i)}(x, y), \ldots, \hat{\phi}_S^{(i)}(x, y) )
    \label{eq:sdf_polygon}
\end{equation}

Recall that the SDF defines the shortest distance of any point in the domain from the boundary of the object. The SDF is negative for any point in the interior of the polygon, zero on the boundary and positive for points exterior of the polygon.
Additionally, note that while the exact SDF satisfies $||\nabla\phi(\bm{x})|| = 1 \; , \; \forall \bm{x}$ it is not necessarily satisfied in the approximate SDF \cite{deng2020cvxnet}. We then compute the density field of the polygon $\hat{\rho}(\cdot)$ from $\phi(\cdot)$ using the Sigmoid function:

\begin{equation}
    \hat{\rho}^{(i)}(x,y; \phi^{(i)}) = \frac{1}{1 + e^{-\beta \phi^{(i)}(x,y)}}
    \label{eq:density_projection}
\end{equation}

\begin{figure}[h]
 	\begin{center}
		\includegraphics[scale=0.5,trim={0 0 0 0},clip]{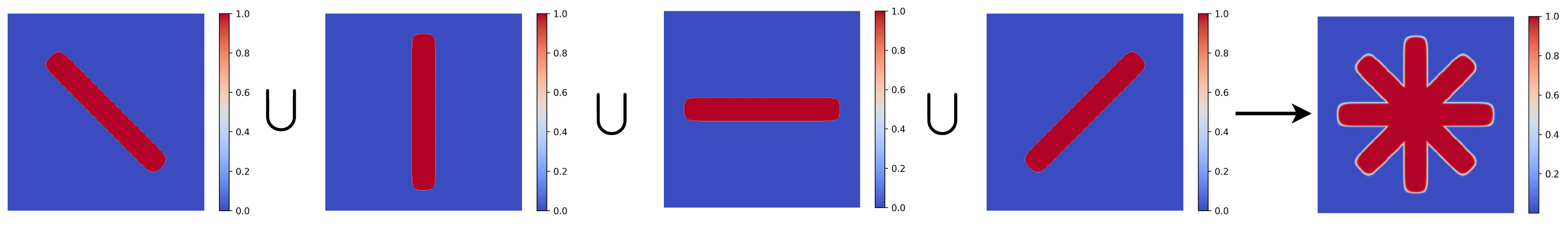}
 		\caption{The density field of the design is computed from the density field of the polygons.}
 		\label{fig:density_union}
	\end{center}
 \end{figure}

 Observe that negative values of the SDF are projected to a density of one and positive values to zero. SDF Values close to zero are projected to intermediate density values $ \rho \in (0,1)$, with the sharpness of transition controlled by the hyperparameter:  $\beta$. Finally, with the density field of each of the polygons $\{\hat{\rho}^{(1)},\ldots,\hat{\rho}^{(K)}\}$ computed, the density field of the design can be derived via a $q^{th}$ norm formulation, as (Eq.~\ref{eq:density_union}, illustrated in Fig.~\ref{fig:density_union}):

\begin{equation}
    \rho(x,y; \hat{\rho}^{(1)}, \ldots, \hat{\rho}^{(K)}) = \bigg(\sum\limits_{k=1}^{K} \big(\hat{\rho}^{(k)}(x,y)\big)^{q}\bigg)^{\frac{1}{q}}
    \label{eq:density_union}
\end{equation}

\subsection{Finite Element Analysis}
\label{sec:method_FEA}

Having obtained the design's density field, we use the Solid Isotropic Material Penalization (SIMP) scheme to determine the Young's modulus $E$. With $\rho_e$ being the density computed at the coordinates of the center of element $e$ for all $N_e$ elements, we obtain the corresponding $E_e$ as:

\begin{equation}
    E_e(\rho_e) = E_{min} + (E_0 - E_{min})\rho_e^p
    \label{eq:SIMP_material}
\end{equation}

Where $E_0$ is the Young's modulus of the base material, $E_{min}$ is a small constant added to prevent a singular global stiffness matrix, and $p$ is the SIMP penalty. 

For finite element analysis, we solve the linear elasticity problem using a structured mesh with bilinear quad element. One can evaluate the element stiffness matrix as

\begin{equation}
[K_e] = E_e \int\limits_{\Omega_e} [B]^T[D_0][B] d \Omega_e
\label{eq:elem_stiffness}
\end{equation}

where $[B]$ is the strain matrix, and $[D_0]$ is the constitutive matrix with a Young’s modulus of unity and an assumed Poisson's ratio (see Table~\ref{table:defaultParameters}). This is followed by the assembly of the global stiffness matrix $\bm{K}$. Finally, with the imposed nodal force vector $\bm{f}$, we solve the state equation for the nodal displacements $\bm{u} = \bm{K}^{-1}\bm{f}$.

\subsection{Optimization Problem}
\label{sec:method_opt}

In this work, we consider the TO of structural components. The various parts of the optimization formulation are discussed below.

\paragraph{Optimizer:} We employ the method of moving asymptotes (MMA) \cite{svanberg1987MMA} to perform the design updates. Specifically, we use a Python implementation \cite{mmapygithub} with all default parameters corresponding to the version of MMA presented in \cite{svanberg2007mma}. For the choices of the MMA move limit, the step tolerance, and the Karush-Kuhn-Tucker (KKT) tolerance \cite{nocedal1999numerical} see Table~\ref{table:defaultParameters}.

\paragraph{Optimization variables:} The design is defined by the polygon's center coordinates, angular offset and plane distances, totalling ${K(S+3)}$ free parameters (Section~\ref{sec:method_designParameterization}). For optimization, we define an augmented normalized design vector $\bm{z} = [\bm{z}_{c_x}, \bm{z}_{c_y}, \bm{z}_{\alpha}, \bm{z}_{d},]$ where the entries are in $[0,1]$. With $ \barbelow{c}_x, \barbelow{c}_y, \barbelow{\alpha}, \barbelow{d} $ being the lower bound and  $ \overline{c}_x, \overline{c}_y, \overline{\alpha}, \overline{d} $ being the upper bound on the x-center, y-center, angular offset and planar distances respectively,  we can retrieve the unnormalized x-center as  $\bm{c_x} \leftarrow \barbelow{c}_x + (\overline{c}_x - \barbelow{c}_x)\bm{z_{c_x}}$ and so on.

\paragraph{Objective:} We adopt a simple compliance minimization problem. With the nodal displacements $\bm{u}$ computed (Section~\ref{sec:method_FEA}), and with $\bm{f}$ being the imposed load, the compliance is  computed as $J = \bm{f}^T \bm{u}$.

\paragraph{Volume constraint:} The designs are subjected to a total volume constraint (Eqs.~\ref{eq:vol_cons}) where $v_f^*$ is the maximum allowed volume fraction, $\rho_e$ are the elemental densities, $v_e$ are the element areas, $\Omega_0$ is the volume of the domain

\begin{equation}
    g_v \equiv \frac{\sum\limits_{e=1}^{N_e} \rho_e v_e}{v_f^* |\Omega_0|} - 1\leq 0
    \label{eq:vol_cons}
\end{equation}

\begin{figure}[h]
 	\begin{center}
		\includegraphics[scale=0.5,trim={0 0 0 0},clip]{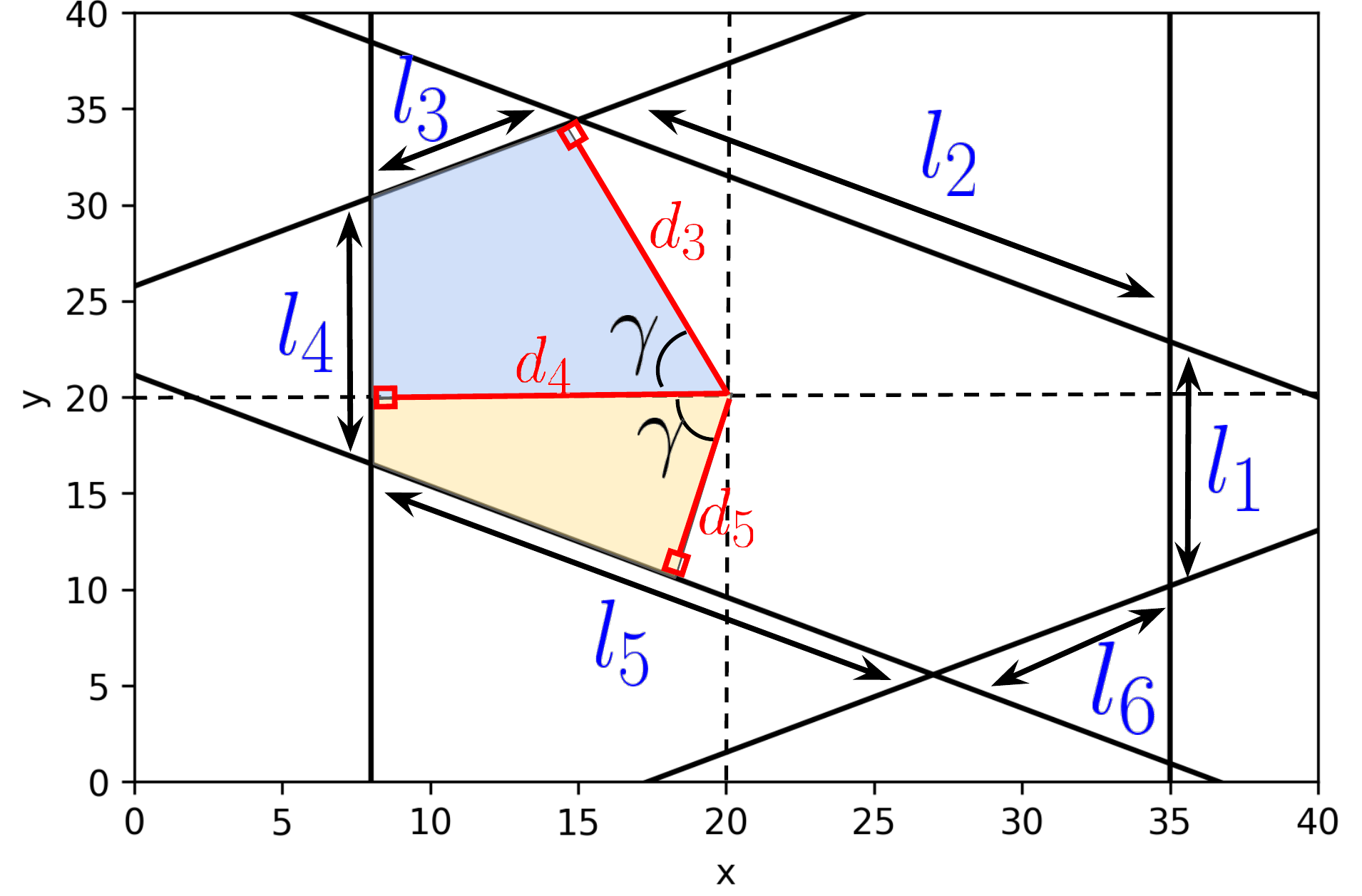}
 		\caption{The length of the polygon edges ($l_1, l_2,\ldots,l_6$) are to be constrained by a minimum bound. The length of any edge (here $l_4$) can be arrived by reasoning from the two encasing quadrilaterals (shown in blue and yellow).}
 		\label{fig:min_feature_size_illustration}
	\end{center}
 \end{figure}

\paragraph{Minimum length constraint:} Optionally, one may wish to impose constraints on the feature sizes of the polygons to promote manufacturability. In particular, we consider a minimum bound imposed on the edge lengths of the polygon. Consider Fig.~\ref{fig:min_feature_size_illustration} illustrating a single hexagon with edge lengths  $l_1, l_2,\ldots l_6$. In general, given our parameterization (Sec.~\ref{sec:method_designParameterization}) these lengths can be derived considering the two enclosing quadrilaterals as (Eq.~\ref{eq:edge_length_polygon}):

\begin{equation}
    l_i^{(j)} = \frac{1}{\sin\gamma}[d_{i+1}^{(j)} + d_{i-1}^{(j)} - 2d_i^{(j)} \cos\gamma]
    \label{eq:edge_length_polygon}
\end{equation}

Where $\gamma = \frac{2\pi}{S}$, and $d_i^{(j)}$, $d_{i-1}^{(j)}$, $d_{i+1}^{(j)}$ are the candidate $i$ and cyclic adjacent perpendicular distances of the $j^{th}$ polygon respectively. With the lengths of all  edges of all the polygons computed, the approximate length of the smallest edge can be expressed differentiably as:

\begin{equation}
    l_{min} = \underset{\forall i, \forall j}{\text{min}} \; l_i^{(j)} \approx -\underset{\forall i, \forall j}{LSE}(-l_i^{(j)})
    \label{eq:smooth_min_length}
\end{equation}

Finally, with $l^*$ being the imposed minimum edge length, one may express the minimum edge length constraint as:

\begin{equation}
    g(l) \equiv 1 - \frac{l_{min}}{l^*} \leq 0
    \label{eq:min_length_scale_cons}
\end{equation}

A similar argument can be used to construct a constraint on the maximum feature size of the polygon, but is not pursued in the current work.

Collecting the objective, state equation, volume, min length and box constraints, the optimization problem can be written as (Eqs.~\ref{eq:optimization_eq}):

\begin{subequations}
	\label{eq:optimization_eq}
	\begin{align}
		& \underset{\bm{z}}{\text{minimize}}
		& &J = \bm{f}^T \bm{u}(\bm{z}) \label{eq:opt_objective}\\
		& \text{subject to}
		& & \bm{K}(\bm{z})\bm{u} = \bm{f}\label{eq:opt_fea}\\
		& & & g_{v} (\bm{z}) \leq 0 \label{eq:opt_volCons}\\
            & & & g_{l} (\bm{z}) \leq 0 \label{eq:opt_minLenCons}\\
            & & &  0 \leq z_i \leq 1 \quad \forall i \label{eq:opt_boxCons}
	\end{align}
\end{subequations}

\subsection{Sensitivity Analysis}
\label{sec:sensitivity_analysis}
A critical ingredient for the update schemes in gradient-based optimization is the sensitivity, i.e., derivative, of the objective and constraint(s) with respect to the optimization parameters. Typically the sensitivity analysis is carried out manually. This can be laborious and error-prone, especially for non-trivial objectives. Here, by expressing all our computations including computing the SDF, projection, FEA, objectives and constraints in JAX \cite{jax2018github}, we use automatic differentiation (AD) capabilities to completely automate this step \cite{chandrasekhar2021auto, ian2020AD}. We can then express the derivative of the objective, volume 
 and minimum length scale constraint with respect to the optimization variables in Eqs.~\ref{eq:sens_objective}, ~\ref{eq:sens_volCons} and ~\ref{eq:sens_MinLenCons} respectively.

\begin{equation}
    \frac{\partial J}{\partial \bm{z}} = \frac{\partial J}{\partial \bm{u}} \frac{\partial \bm{u}}{\partial \bm{K}} \frac{\partial \bm{K}}{\partial \bm{E}}\frac{\partial \bm{E}}{\partial \bm{\rho}} \frac{\partial \bm{\rho}}{\partial \bm{\hat{\rho}}} \frac{\partial \bm{\hat{\rho}}}{\partial \bm{\phi}} \frac{\partial \bm{\phi}}{\partial \bm{\hat{\phi}}} \frac{\partial \bm{\hat{\phi}}}{\partial \bm{p}}\frac{\partial \bm{p}}{\partial \bm{z}}
    \label{eq:sens_objective}
\end{equation}

\begin{equation}
    \frac{\partial g_v}{\partial \bm{z}} = \frac{\partial g_v}{\partial \bm{\rho}} \frac{\partial \bm{\rho}}{\partial \bm{\hat{\rho}}} \frac{\partial \bm{\hat{\rho}}}{\partial \bm{\phi}} \frac{\partial \bm{\phi}}{\partial \bm{\hat{\phi}}} \frac{\partial \bm{\hat{\phi}}}{\partial \bm{p}}\frac{\partial \bm{p}}{\partial \bm{z}}
    \label{eq:sens_volCons}
\end{equation}

\begin{equation}
    \frac{\partial g_l}{\partial \bm{z}} = \frac{\partial g_l}{\partial \bm{p}}\frac{\partial \bm{p}}{\partial \bm{z}}
    \label{eq:sens_MinLenCons}
\end{equation}

\subsection{Algorithm}
\label{sec:method_algo}

\begin{figure}[h]
 	\begin{center}
		\includegraphics[scale=0.6,trim={0 0 0 0},clip]{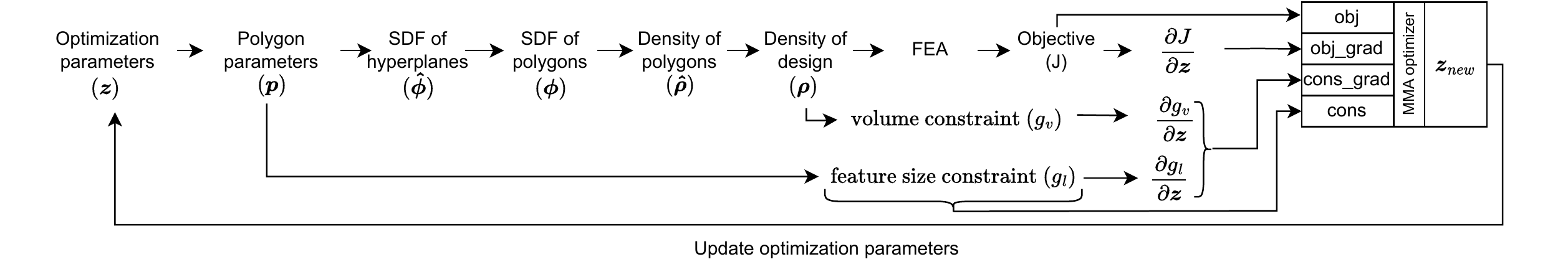}
 		\caption{Optimization loop of the proposed framework.}
 		\label{fig:flowchart}
	\end{center}
 \end{figure}

The complete algorithm of the proposed framework is summarized in Algorithm 1, and schematically depicted in Fig.~\ref{fig:flowchart}.

\begin{algorithm}
    \caption{PolyTO: Topology Optimization using Convex Polygons}\label{alg:PolyTO}
    \begin{algorithmic}
    \State $i \leftarrow 0$ \Comment{Iteration counter}
    \State $\bm{z}^{(0)} \leftarrow \bm{p}^{(0)}$ \Comment{Initialize opt. variables}
    \State $\text{converged = False}$
    \While{$\text{converged == False}$}
    \State $\bm{p}^{(i)} \leftarrow \bm{z}^{(i)}$ \Comment{Get design from opt. variables}
    \State {Compute minimum length constraint $g_l$ from $\bm{p}^{(i)}$} \Comment{Eqs.~\ref{eq:edge_length_polygon}, ~\ref{eq:smooth_min_length}, ~\ref{eq:min_length_scale_cons}}
    \For{$k = 1,\ldots,K$} \Comment{every polygon}
    \For{$s = 1,\ldots,S$} \Comment{every half-space}
    \For{$e=1,\ldots,N_e$} \Comment{every element}
    \State {Compute signed distance $\hat{\phi}^{(k)}_{s}(\bm{x}_e)$ from $\bm{x}_e$ to side $s$ of polygon $k$} \Comment{Eq.~\ref{eq:sdf_hyperplane}}
    \EndFor
    \EndFor
    \State {Compute signed distance of polygon  $\phi^{(k)}(\bm{x}_e)$} \Comment{Eq.~\ref{eq:sdf_polygon}}
    \State {Compute projected density $\hat{\rho}^{(k)}(\bm{x}_e)$ of polygons} \Comment{Eq.~\ref{eq:density_projection}}
    
    \EndFor

    \For{$e=1,\ldots,N_e$}
    \State {Compute union of density $\rho_e$}
    \Comment{Eq.~\ref{eq:density_union}}
    \State {Compute Modulus of Elasticity $E_e$}\Comment{Eq.~\ref{eq:SIMP_material}}
    \State {Compute element stiffness matrix $\bm{K}_e$}\Comment{Eq.~\ref{eq:elem_stiffness}}
    \EndFor
    \State {Compute assembled stiffness matrix $\bm{K}$ and solve $\bm{u} = \bm{K}^{-1}\bm{f} $ }\Comment{Section~\ref{sec:method_FEA}}
    \State{Compute Objective $J$ and volume constraint $g_v$} \Comment{Section~\ref{sec:method_opt}}
    \State{Compute $\nabla_{\bm{z}}J$, $\nabla_{\bm{z}}g_v$, $\nabla_{\bm{z}}g_l$ } \Comment{Eqs.~\ref{eq:sens_objective}, ~\ref{eq:sens_volCons}, ~\ref{eq:sens_MinLenCons}}
    \State{$\bm{z^{(i+1)}} \leftarrow \text{MMA}(J, g_v, g_l, \nabla_{\bm{z}}J, \nabla_{\bm{z}}g_v, \nabla_{\bm{z}}g_l)$} \Comment{Update design}
    \State{ $i \leftarrow i + 1 $} \Comment{Increment counter}
    \State{$\Delta_{s} \leftarrow ||\bm{z}^{(i)} - \bm{z}^{(i-1)}||$} \Comment{Compute step tolerance}
    \State {Compute norm of KKT condition $\Delta_{kkt}$}
    \If{$k \geq \text{maxIter}$  or $\Delta_{s} \leq \text{stepTol}$ or $\Delta_{kkt} \leq \text{kktTol}$} \Comment{Convergence criteria}
    \State {converged = True}
    \EndIf
    \EndWhile

    \end{algorithmic}
\end{algorithm}

%% file: results.tex
In this section, we conduct several experiments to illustrate the proposed framework. The default parameters in the experiments are summarized in Table~\ref{table:defaultParameters}. The minimum length scale constraint is not utilized except for experiments in Sec.~\ref{sec:result_minLength}. All experiments are conducted on a MacBook M1 Pro, using the JAX library \cite{jax2018github} in Python. 

\begin{table}[!t]
	\caption{Default parameters}
	\begin{center}
		\begin{tabular}{  r | p{62mm}  }
			
        Parameter & Description and default value \\ \hline
        $E$, $\nu$ & Isotropic material with Young's Modulus $E = 1$ and Poisson's ratio $\nu = 0.3$ \\
   
        nelx, nely & Number of FEA mesh elements of (100, 50)  along $x$ and $y$ respectively \\

        $l_x$, $l_y$ & Length of domain (60, 30) along $x$ and $y$ respectively \\

        $(\underbar{c}_x$, $\overline{c}_x)$ & Polygon x-center range $(0, l_x)$ \\

        $(\underbar{c}_y$, $\overline{c}_y)$ & Polygon y-center range $(0, l_y)$ \\

        $\underbar{d}$, $\overline{d}$ & Polygon face offset range $(0, 0.5l_x)$ \\

        $p$ & SIMP penalty parameter = 3 \\

        $q$ & Density union penalty parameter = 8 \\

        $\beta$ & Sharpness param in density projection = 10 \\
        maxIter & Max number of optim iters = $250$\\
        
        move limit & MMA  step size = $5E-2$\\
   
        kkt tolerance & MMA convergence criteria = $1E-6$ \\
                
        step tolerance & MMA convergence criteria = $1E-6$ \\
			
		\end{tabular}
	\end{center}
	
	\label{table:defaultParameters}
\end{table}

\subsection{Validation}
\label{sec:result_validation}

\begin{figure}[htbp]
 	\begin{center}
		\includegraphics[scale=2.,trim={0 0 0 0},clip]{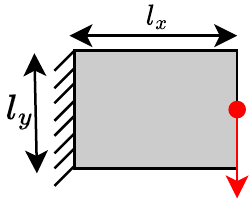}
 		\caption{Mid-cantilever beam}
 		\label{fig:mid_cant}
	\end{center}
 \end{figure}

\begin{figure}[htbp]
 	\begin{center}
		\includegraphics[scale=0.35,trim={0 0 0 0},clip]{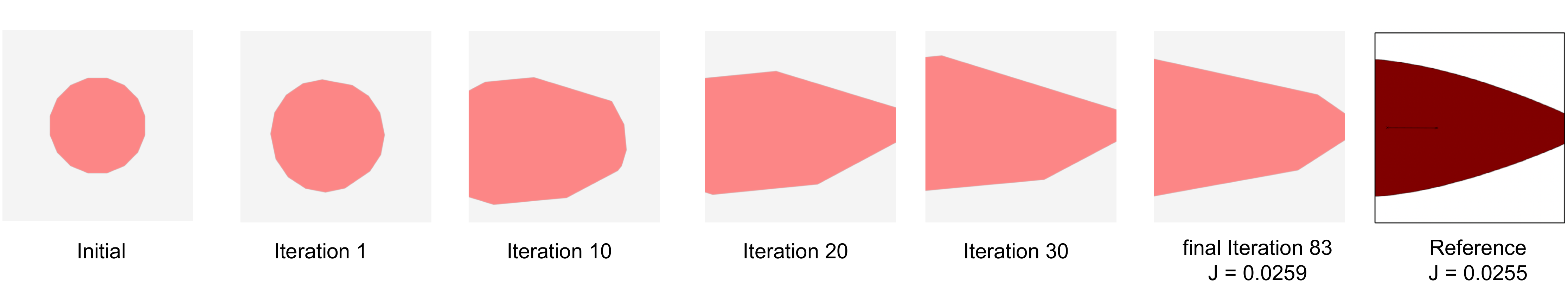}
 		\caption{Design evolution for a mid-cantilever beam with a single polygon. The final result is compared with the result in Ref.~\cite{norato2018Supershapes}}
 	\label{fig:single_blob_evolution_mid_cant}
	\end{center}
 \end{figure}

We consider the example presented in \cite{norato2018Supershapes} where the design of a mid-cantilever beam (Fig.~\ref{fig:mid_cant}) with $l_x = l_y = 2$ for $v_f^* = 0.5$ was considered. The design is parameterized with one polygon with $S=16$. We report the evolution and the final topology in Fig.~\ref{fig:single_blob_evolution_mid_cant}. The topology from the reference is also illustrated. We observe that the optimized topology and the final compliance agrees with the reference, validating our method with established work. Note that as in \cite{norato2018Supershapes}, only the portion of the polygon within the design domain is reported as the design outside the domain has no effect on the objective and constraints. Constraining the polygon to lie entirely within the design domain \cite{kang2016structural} will be considered in the future.

\begin{figure}[htbp]
 	\begin{center}
		\includegraphics[scale=0.5,trim={0 0 0 0},clip]{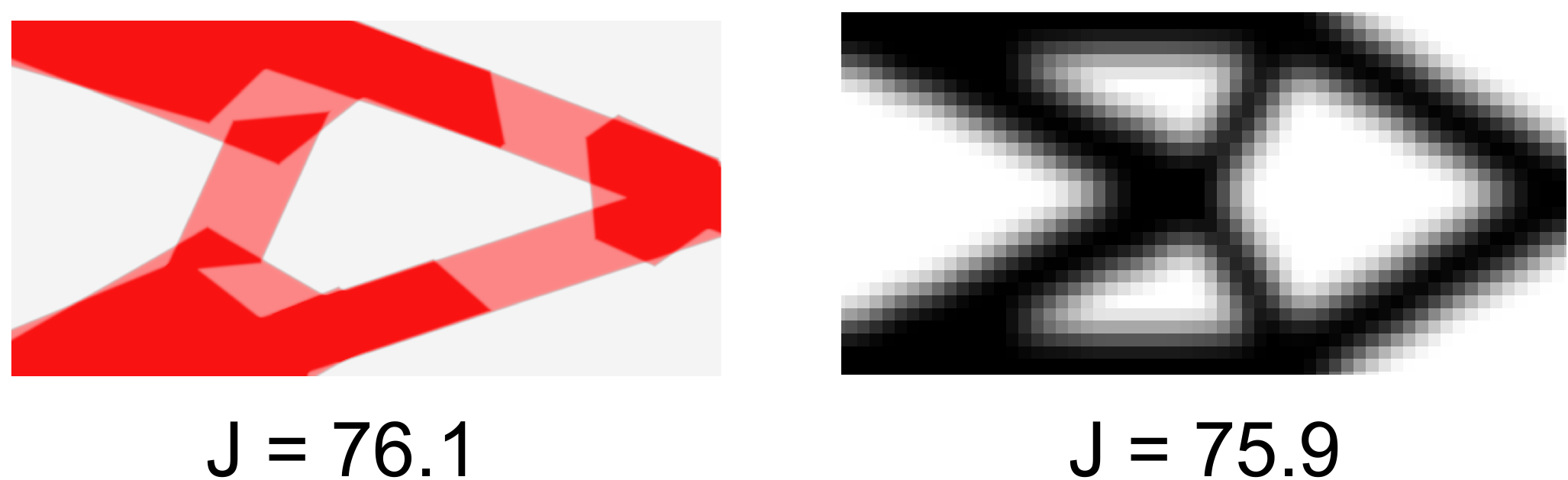}
 		\caption{Comparison of topology and compliance with \cite{sigmund200199}.}
 		\label{fig:validation_top99}
	\end{center}
 \end{figure}

 Next, we compare our result with that obtained from a density based method. Fig.~\ref{fig:validation_top99} illustrates the design obtained using the proposed framework ($K=12$ and $S=6$) with that obtained using \cite{andreassen2011Top88} for a mid-cantilever beam with the default parameters and $v_f^* = 0.5$. We observe that the optimized designs have similar performance. Note that the transparency of the polygons are purely illustrative. 
 
\begin{figure}[htbp]
 	\begin{center}
		\includegraphics[scale=0.5,trim={0 0 0 0},clip]{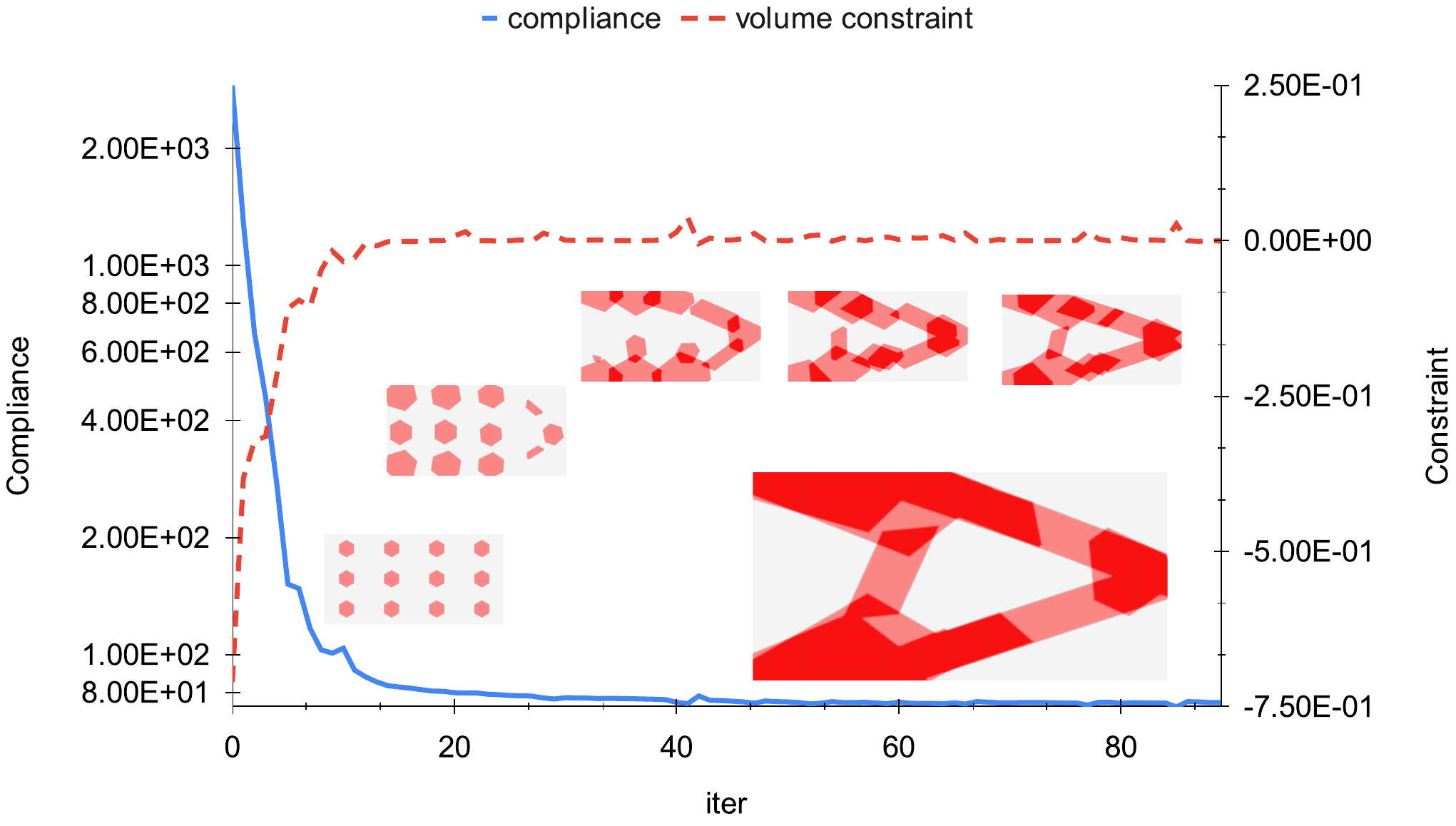}
 		\caption{Convergence of compliance (objective) and volume constraint for a mid-cantilever beam. Designs at the 0th, 5th, 10th, 30th, 50th and final iteration are illustrated.}
 		\label{fig:convergence}
	\end{center}
 \end{figure}

 We also illustrate the convergence of compliance, and volume constraint for the problem in Fig.~\ref{fig:convergence}. Snapshots of the design at various instances of the optimization are also showcased. Similar convergence behavior was observed for all other examples.

\subsection{Pareto trade-off}
 \label{sec:result_pareto}
 
\begin{figure}[htbp]
 	\begin{center}
		\includegraphics[scale=1.5,trim={0 0 0 0},clip]{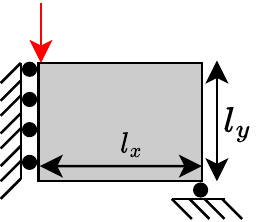}
 		\caption{MBB beam.}
 		\label{fig:MBB}
	\end{center}
 \end{figure}

\begin{figure}[htbp]
 	\begin{center}
		\includegraphics[scale=0.25,trim={0 0 0 0},clip]{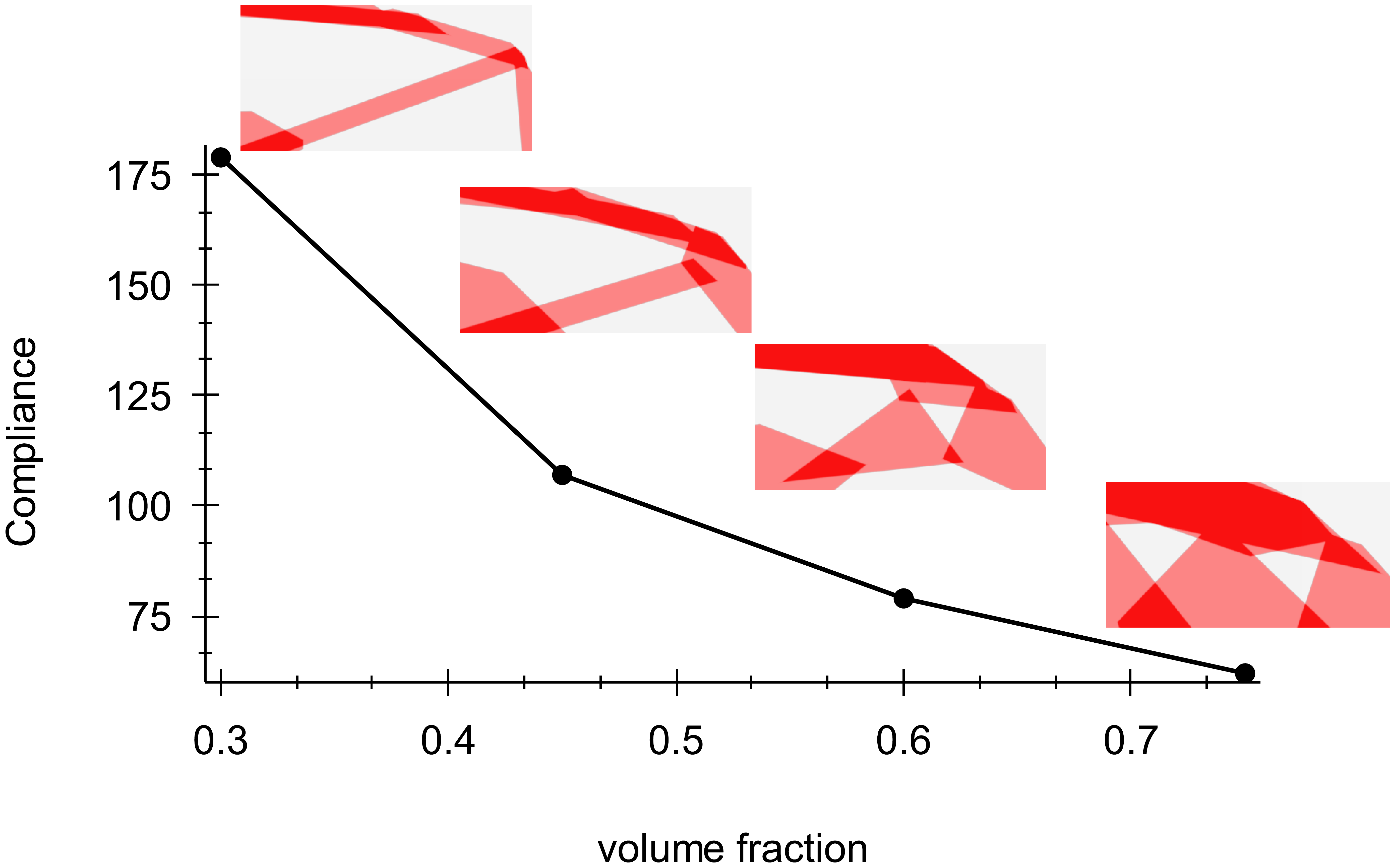}
 		\caption{Trade-off between compliance and volume.}
 		\label{fig:pareto}
	\end{center}
 \end{figure}

An important consideration during the design stage is in exploring the Pareto-front, evaluating the trade-offs for various design choices. Consider the MBB beam problem in Fig.~\ref{fig:MBB}. With $K=6$ and $S=12$, we explore the trade-off between the structure's compliance and volume fraction. Fig.~\ref{fig:pareto} illustrates the obtained pareto-optimal front and topologies at various instances. As expected, we observe increasing compliance as we lower the allowed volume fraction.

\subsection{Effect of initialization}
\label{sec:result_initialization}

\begin{figure}
     \centering
     \begin{subfigure}[b]{\textwidth}
         \centering
         \includegraphics[scale=0.5,trim={0 0 0 0},clip]{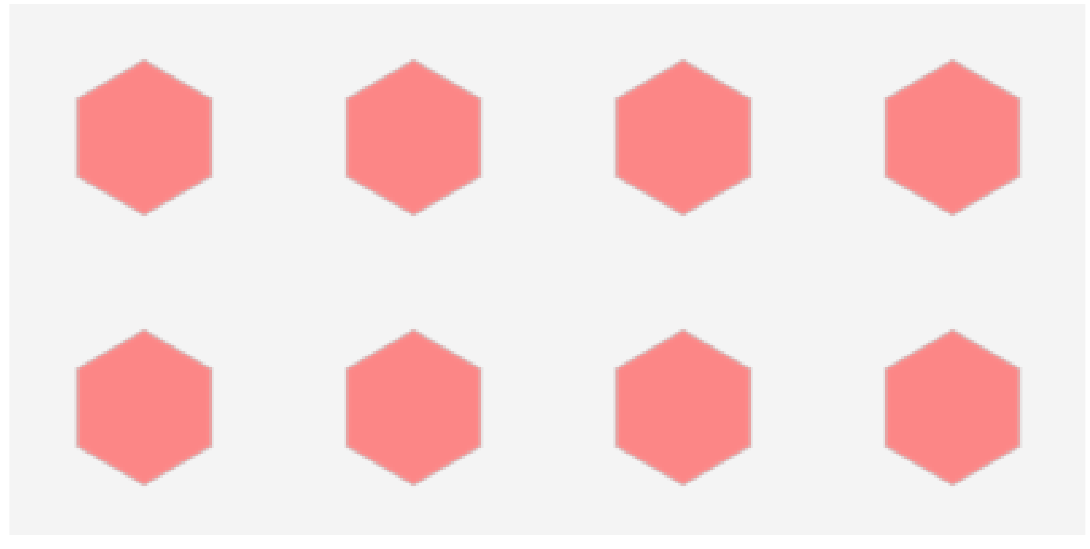}
         \caption{Initial design consisting of $(K=)$ 8 hexagons $(S=6)$ on a $4\times 2$ grid.}
         \label{fig:polygon_init_grid}
     \end{subfigure}
     \vfill
     \begin{subfigure}[b]{\textwidth}
         \centering
         \includegraphics[scale=0.5,trim={0 0 0 0},clip]{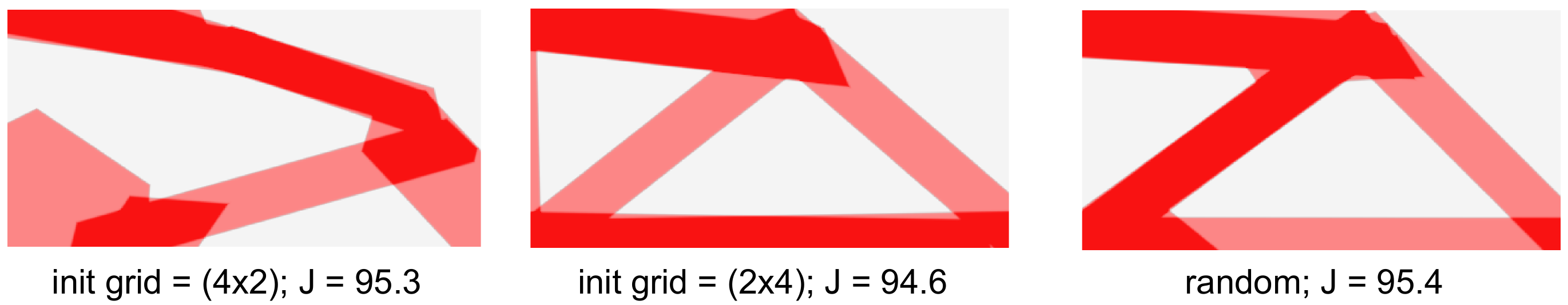}
         \caption{MBB beam design with different initialization.}
         \label{fig:initialization_MBB}
     \end{subfigure}
        \caption{Effect of initialization on the resulting design and compliance.}
        \label{fig:initialization_effect}
\end{figure}

In this experiment, we consider the effect of the initial design on the optimized result. Owing to the non-convex, non-linear nature of the our optimization problem, we expect the different initialization to result in different in local optima. We typically start the optimization with a grid of equi-spaced, equi-sized polygons. For instance, consider Fig.~\ref{fig:polygon_init_grid} where 8 hexagons are initialized on a $4 \times 2$ grid. In this example, we consider the design of a MBB beam (Fig.~\ref{fig:MBB}) with hyper-parameters $K=8$, $S=6$ and $v_f^* = 0.5$. We compare the resulting topology and performance for different initial configurations in Fig.~\ref{fig:initialization_MBB}. We observe that while obtaining varied topologies, the performances are similar.

 \subsection{Varying polygon parameters}
 \label{sec:result_polygonParam}

Next, we consider the effect of two hyper-parameters, namely the number of polygons $(K)$ and number of hyper-planes in a polygon $(S)$. One may control the complexity design by varying the number of polygons and the number half-spaces that parameterize a polygon. Consider once again the MBB beam as shown in Fig.~\ref{fig:MBB}. With $v_f^* = 0.5$, we consider the effect of varying $K$ and $S$ in Fig.~\ref{fig:polygon_params}. While we observed no considerable differences in performance or topology for our problem in consideration and in general for the compliance minimization of structural members, these parameters may prove important for other objectives and physics.

  \begin{figure}[htbp]
 	\begin{center}
		\includegraphics[scale=0.5,trim={0 0 0 0},clip]{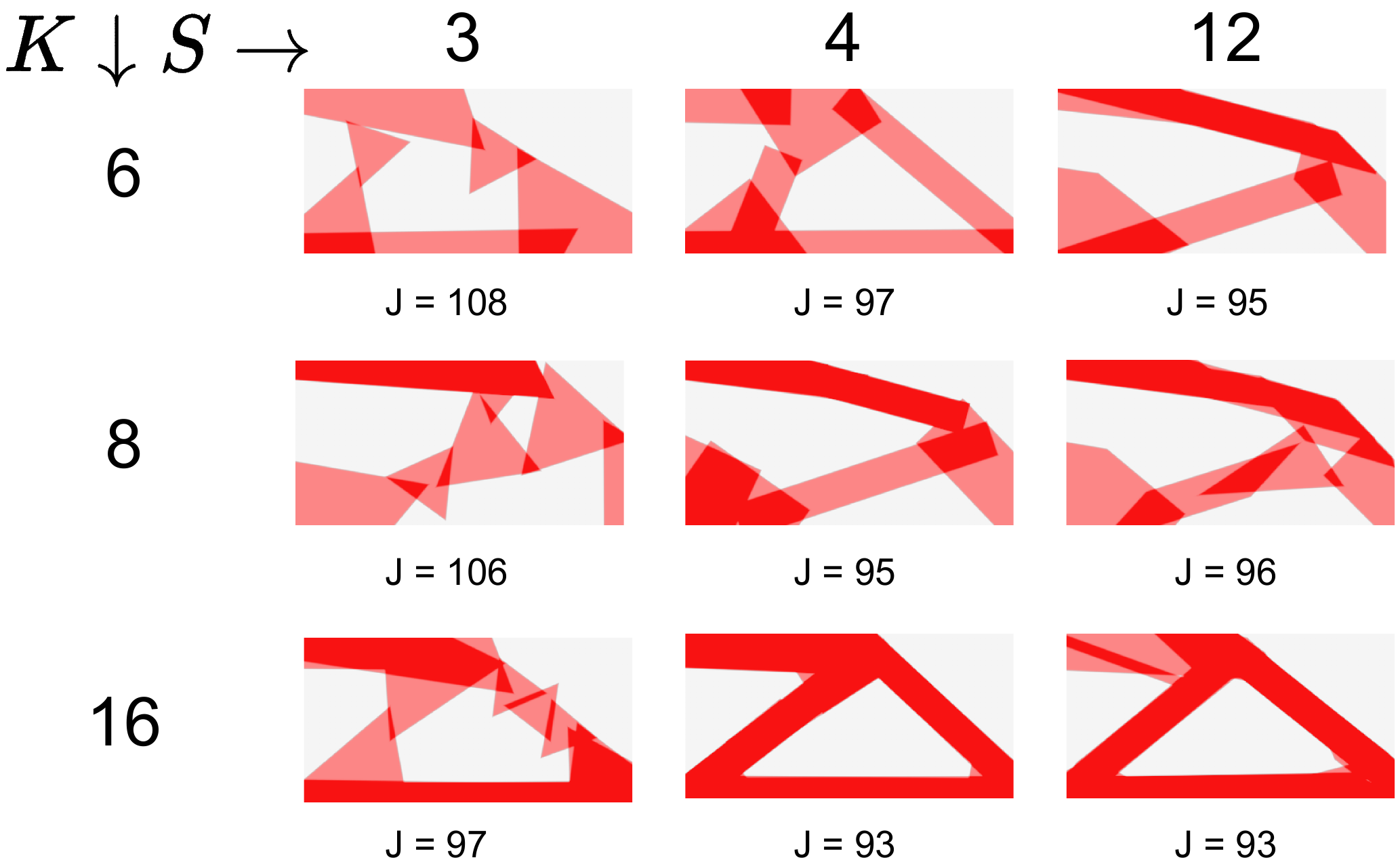}
 		\caption{Effect of number of polygons $(K)$ and number of sides $(S)$ on the resulting topology and performance.}
\label{fig:numPoly_numSides_exploration}
	\end{center}
 \end{figure}

 \subsection{Minimum length scale constraint}
 \label{sec:result_minLength}

\begin{figure}[]
  \begin{subfigure}[b]{\linewidth}
  \centering
    \includegraphics[width=0.5\linewidth]{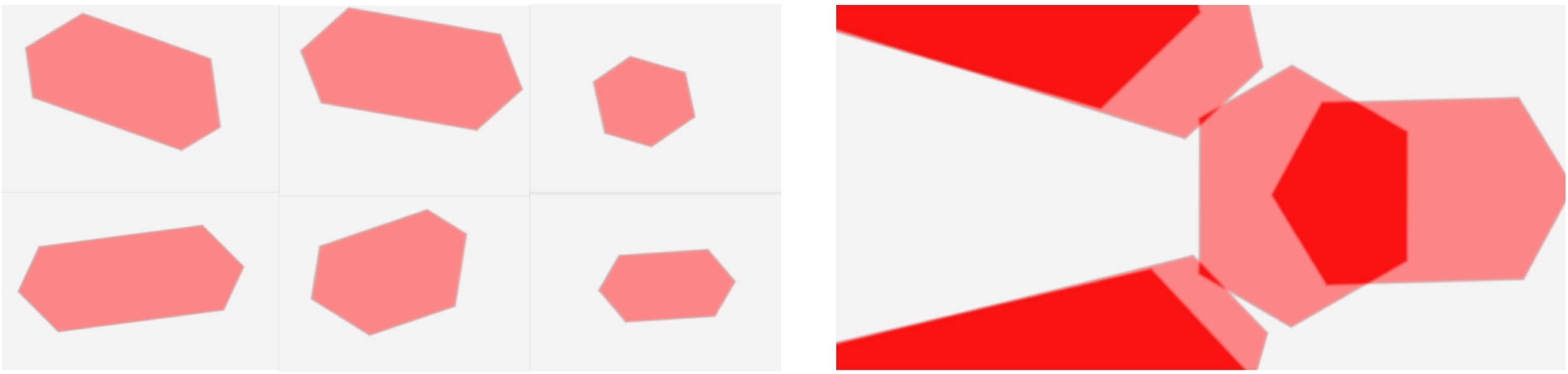}
    \caption{Allowed minimum length $l^* = 8$}
    \label{fig:size_8}
  \end{subfigure}

  \begin{subfigure}[b]{\linewidth}
  \centering
    \includegraphics[width=0.5\linewidth]{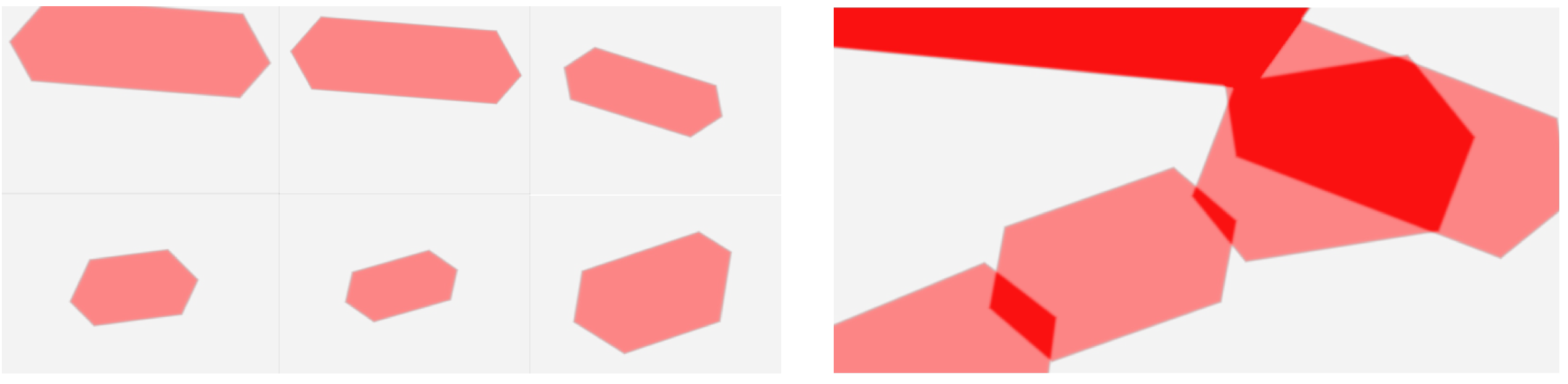}
    \caption{Allowed minimum length of $l^* = 6$}
    \label{fig:size_6}
  \end{subfigure}
  
    \begin{subfigure}[b]{\linewidth}
  \centering
    \includegraphics[width=0.5\linewidth]{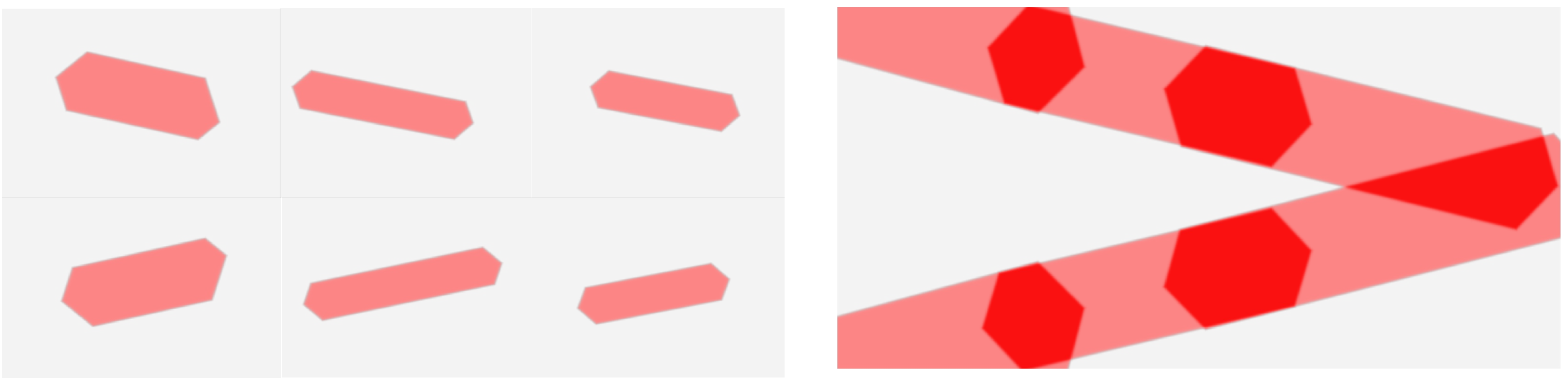}
    \caption{Allowed minimum length of $l^* = 4$}
    \label{fig:size_4}
  \end{subfigure}
  
    \begin{subfigure}[b]{\linewidth}
  \centering
    \includegraphics[width=0.5\linewidth]{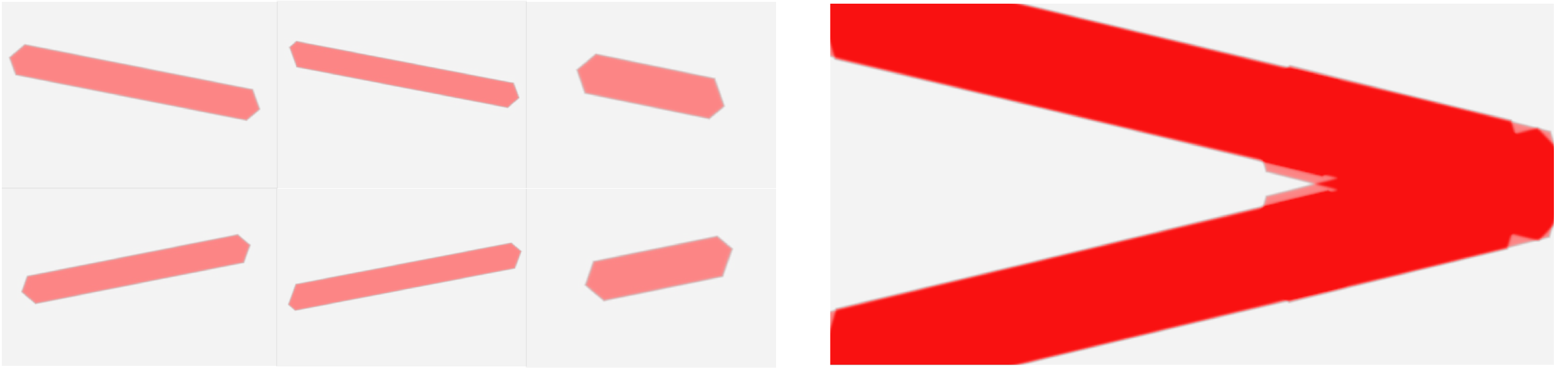}
    \caption{Allowed minimum length of $l^* = 2$}
    \label{fig:size_2}
  \end{subfigure}
  
  \caption{Polygons (left) and resulting design (right) for varying allowed minimum length}
  \label{fig:topology_varying_min_length}
\end{figure}

\begin{figure}[htbp]
 	\begin{center}
		\includegraphics[scale=0.5,trim={0 0 0 0},clip]{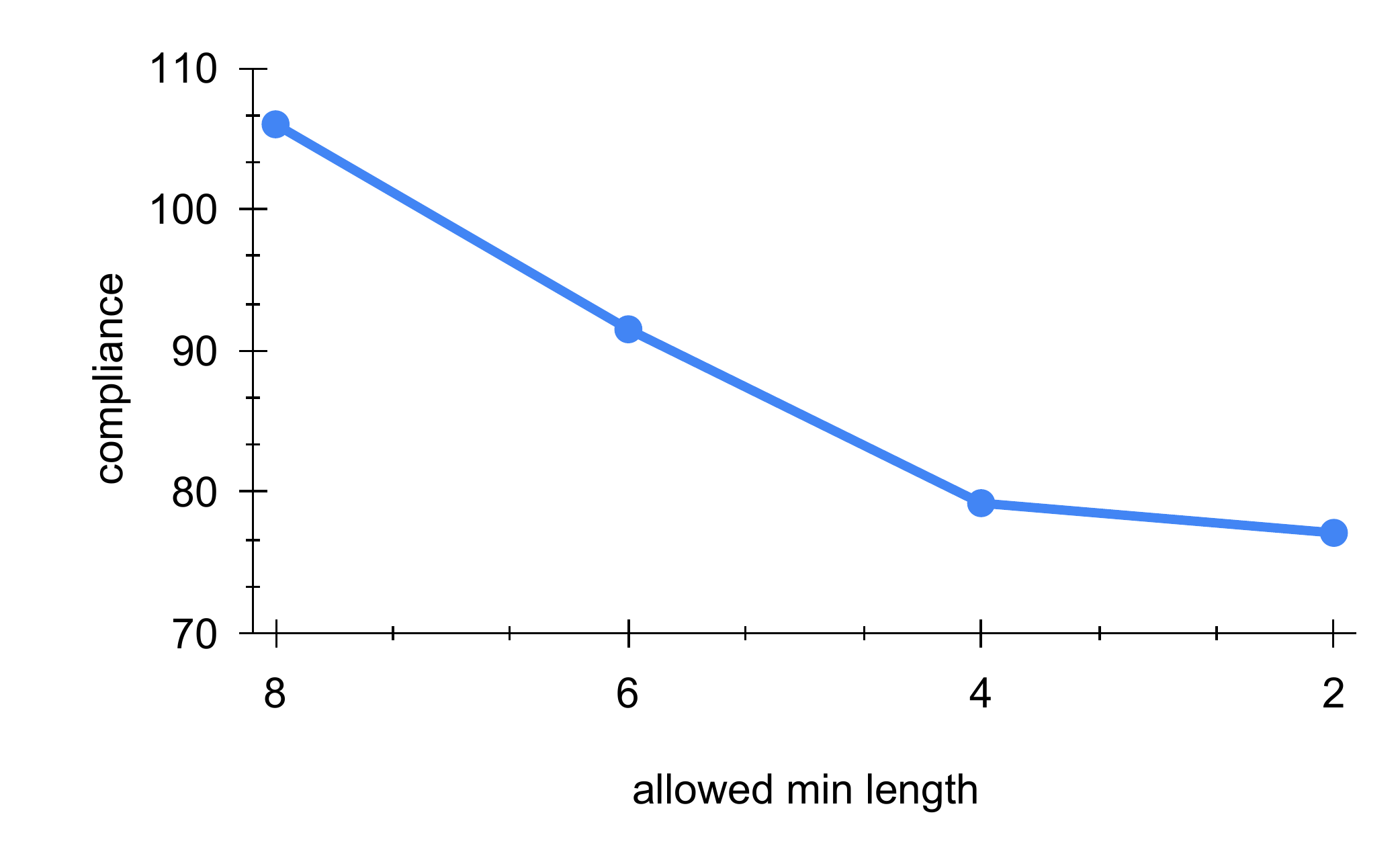}
 		\caption{Effect of allowed minimum length on the performance.}
\label{fig:min_length_exploration}
	\end{center}
 \end{figure}

Next, we consider the effect of the minimum length scale constraint (Eq.~\ref{eq:min_length_scale_cons}) on the resulting topology and performance. Consider once again the mid-cantilever 
 in Fig.~\ref{fig:mid_cant}. Fig.~\ref{fig:topology_varying_min_length} showcases the topology along with the populating polygons for varying allowed minimum length $l^*$ with $K=S=6$ and $v_f^* = 0.5$. Note that a larger allowed minimum length corresponds to a more stricter constraint. As expected, we notice larger polygons with a loss in performance (Fig.~\ref{fig:min_length_exploration}) as we enforce stricter constraints. 

%% file: conclusion.tex
This paper demonstrates a method to perform TO where the design is parameterized by convex polygons. A natural extension to feature mapping methods in TO, polygons accommodates more varied shapes than simpler primitives like bars and plates. This allows for increased design freedom, leading to more optimal solutions while still retaining the advantages of feature mapping methods. The complexity of the geometry can be controlled by the number of polygons and number of half-spaces parameterizing each polygon. We presented the framework in the context of structural optimization to obtain stiff designs subject to a volume constraint. The feature sizes were also derived and constrained analytically. The numerical examples show the method produces designs with performances competitive to those obtained with pixel based parameterization.

We assumed the number of sides and number of polygons as hyper-parameters given a-priori to the optimization. Incorporating these as part of the optimization parameters may lead to better and more automated designs. Further, to ensure we always obtain a bounded and non-empty polygons, we fixed the inter angles between the half-spaces. Future work will explore techniques for free parameterization of the angles. While we limited our discussion to 2D, the technique can be extended to 3D trivially. While we explored feature size control of individual primitives, overlapping primitives can lead to features larger than the specified minimum size. This can be overcome by imposing overlap constraints in addition to the minimum feature size constraint \cite{kang2013integrated}. To further facilitate manufacturability, other constraints such as overhang constraint \cite{qian2017OverhangCons} etc will be explored in the future. Further, the polygons can be eroded and dilated by simply changing the half-space's distances. This can be then be incorporated in a robust analysis formulation \cite{wang2011RobustTO}. Finally, we note that computing the SDF can be expensive. This may be overcome by adopting machine learning based approaches \cite{park2019deepsdf} to obtain the SDF.